\documentclass[11pt]{article}
\usepackage{amssymb}

\def\sgn{\hbox{sgn}}

\newtheorem{theorem}{Theorem}[section]

\newtheorem{lem}{Lemma}[section]
\newtheorem{cor}{Corollary}[section]

\newtheorem{rem}{Remark}[section]
\newcommand{\R}{{ I\!\!R}}

\newcommand{\ep}{\varepsilon}

\newcommand{\re}[1]{(\ref{#1})}

\begin{document}
\begin{center}
\noindent {\Large \bf{Stability of multipeakons}}
\end{center}
\vskip0.2cm
\begin{center}
\noindent
{\bf Khaled El Dika and  Luc Molinet}\\
{\small
L.A.G.A., Institut Galil\'ee, Universit\'e Paris-Nord,\\
93430 Villetaneuse, France.} \vskip0.3cm
 \noindent
% E-mail :
khaled@math.u-paris13.fr\\
% E-mail :
molinet@math.u-paris13.fr
\end{center}
\vskip0.5cm \noindent {\bf Abstract.} { The Camassa-Holm equation
possesses well-known peaked solitary waves that are called
 peakons.
 Their orbital stability  has been established by
 Constantin and Strauss in \cite{CS1}. We prove here the stability of
 ordered  trains of peakons. We also establish a result on the stability of multipeakons.}

%%%%%%%% INTRO

\section{Introduction}
The  Camassa-Holm equation (C-H)$_\kappa $, $ \kappa\ge 0$,
\begin{equation}
u_t -u_{txx}=-2\kappa u_x - 3 u u_x +2 u_x u_{xx} + u u_{xxx},
\quad (t,x)\in\R^2, \label{CH-kappa}
\end{equation}
can be derived as a model
for the propagation of unidirectional shalow water waves over a
flat bottom  by  writing the Green-Naghdi equations in
Lie-Poisson Hamiltonian form and then making an asymptotic
expansion which keeps the Hamiltonian structure (\cite{CH1},
\cite{Johnson}).  It was also found independently by
Dai \cite{dai} as a model for nonlinear waves in cylindrical
hyperelastic rods and was, in fact, first discovered by the
method of recursive operator by Fokas and Fuchsteiner \cite{FF}
as an example of bi-Hamiltonian equation.

(C-H)$_\kappa $ is completely integrable (see
\cite{CH1},\cite{CH2}). It possesses among  others the following
invariants
\begin{equation}
 E(v)=\int_{\R} v^2(x)+v^2_x(x)
\, dx \mbox{ and } F(v)=\int_{\R} v^3(x)+v(x)v^2_x(x)+2\kappa
v^2(x) \, dx\;  \label{E}
\end{equation}
and can be written in Hamiltonian form as
\begin{equation}
\partial_t E'(u) =-\partial_x F'(u) \; .
\end{equation}
For $ \kappa>0 $ it possesses smooth positive solitary waves $
\varphi_{\kappa,c} $ with speed $ c>2\kappa $, their orbital
stability has been proved in \cite{CS2} by applying the classical
spectral method initiated by Benjamin \cite{Benjamin} (see also
\cite{GSS}). In \cite{EL}, following the general method developed
in  \cite{MMT} (see also \cite{EM}), the authors proved the
stability of ordered trains of such solitary waves. It is worth
recalling that this general method requires principally two
ingredients :  A property of almost monotonicity which says that
for a solution close to $\varphi_{\kappa,c}$, the part of the
energy traveling at the right of
 $\varphi_{\kappa,c}(\cdot -ct) $
is almost time decreasing;
 A dynamical proof of the stability of the solitary
wave using the spectral approach
 (as in \cite{Benjamin} or \cite{GSS} for instance). \vspace*{2mm} \\
In this paper we consider the Camassa-Holm equation in the case $\kappa=0$, that is
\begin{equation}
u_t -u_{txx}= - 3 u u_x +2 u_x u_{xx} + u u_{xxx} , \quad (t,x)\in\R^2 .\label{CH}
\end{equation}
Henceforth, we refer to \re{CH} as the Camassa-Holm equation (C-H).
 \re{CH} possesses also solitary waves but they are non
 smooth and are called peakons.
They  are given by
$$ u(t,x)=\varphi_c(x-ct)=c\varphi(x-ct)=ce^{|x-ct|},\ \ \ c\in\R.$$
Their stability seems not to enter the general framework
mentioned above (see the beginning of Section \ref{section1peak} for
further commentaries on this aspect). However, Constantin and
Strauss \cite{CS1} succeeded in proving their orbital stability by
a direct approach. In this work, following the general strategy
initiated in \cite{MMT}(note that due to  the reasons mentioned
above, the general method of \cite{MMT} is not directly
applicable here ), we combine the monotonicity result proved in
\cite{EM} with localized versions of the estimates established in
\cite{CS1} to derive the stability of the trains of peakons.

Before stating the main result we have to introduce the function space where will
  live our class of solutions to the equation. For $ I $  a finite or infinite interval of $ \R $,
  we denote by $ Y(I)  $ the function space\footnote{$ W^{1,1}(\R) $ is the space of $ L^1(\R) $ functions with derivatives in $ L^1(\R) $ and $ BV(\R) $ is the space of function with bounded variation}
  \begin{equation}
Y(I):= \Bigl\{ u\in C(I;H^1(\R)) \cap L^\infty(I;W^{1,1}(\R)), \;
u_x\in L^\infty(I; BV(\R))\Bigr\} \, . \label{theoweak}
\end{equation}
We are now ready to state our main result.
\begin{theorem} \label{mult-peaks}
Let be given $ N $ velocities $c_1,.., c_N $ such that $0<c_1<c_2<..<c_N $.
There exist  $ \gamma_0 $, $ A>0 $, $ L_0>0 $
 and $ \varepsilon_0>0 $ such that if $ u \in Y([0,T[)$, with $ 0<T\le \infty $,
   is a solution of (C-H) satisfying
 \begin{equation}
 \|u_0-\sum_{j=1}^N \varphi_{c_j}(\cdot-z_j^0) \|_{H^1} \le \varepsilon^2 \label{ini}
 \end{equation}
 for some  $ 0<\varepsilon<\varepsilon_0$ and $ z_j^0-z_{j-1}^0\ge L$,
with $ L>L_0 $, then there exist $x_1(t), ..,x_N(t) $ such that
\begin{equation}
\sup_{[0,T[} \|u(t,\cdot)-\sum_{j=1}^N \varphi_{c_j}(\cdot-x_j(t)) \|_{H^1} \le
A(\sqrt{\varepsilon}+L^{-{1/8}})\;  \label{ini2}
\end{equation}
and
\begin{equation}
x_j(t)-x_{j-1}(t)> L/2, \quad \forall t\in[0,T[ \; .
\end{equation}
 \end{theorem}

As discovered by Camassa and Holm \cite{CH1}, (C-H) possesses also special
solutions called multipeakons given by
$$
u(t,x)=\sum_{i=1}^N p_j(t) e^{-|x-q_j(t)|} \,,
$$
where $ (p_j(t),q_j(t)) $ satisfy the differential system \re{systH}.
 In \cite{Beals} (see also \cite{CH1}), the  asymptotic behavior of the   multipeakons is
  studied. In particular, the limits as $ t $ tends to $ +\infty $ and $ -\infty $
   of $ p_i(t) $ and $ \dot{q_i}(t) $ are determined.
 Combining these asymptotics  with the  preceding theorem we get the following result
on the stability of the variety $ {\cal N} $ of $ H^1(\R) $
defined by
$$
{\cal N}:= \Bigl\{ v=\sum_{i=1}^N p_j e^{-|\cdot-q_j|}, \,
(p_1,..,p_N)\in (\R_+)^N , \, q_1<q_2<..<q_N  \Bigr\} \; .
$$
\begin{cor} \label{cor-mult-peaks}
Let be given $ N $  positive real numbers $ p_1^0,.., p_N^0 $ and $ N
$ real numbers $ q_1^0< ..< q_N^0 $. For any $ B> 0 $ and any  $
\gamma >0 $ there exists $ \alpha>0 $ such that if $ u_0\in
H^1(\R) $ satisfies\footnote{ ${\mathcal M}(\R)$ is the space of  Radon measures on $ \R $ with bounded total variation and  $ {\mathcal M}_+(\R) $ is the subset of non-negative
 measures} $ m_0:=u_0-u_{0,xx} \in {\mathcal M}_+(\R) $ with
\begin{equation}
\|m_0\|_{\cal M}\le B \quad  \mbox{ and }\quad
\|u_0-\sum_{j=1}^N p_j^0 \exp (\cdot-q_j^0) \|_{H^1}\le \alpha
\label{ini3}
\end{equation}
 then
\begin{equation}
\forall t\in\R, \quad \inf_{P\in (\R_+)^N,Q\in \R^N}
\|u(t,\cdot)-\sum_{j=1}^N p_j \exp (\cdot-q_j) \|_{H^1} \le
\gamma\; . \label{ini33}
\end{equation}
Moreover, there exists $ T>0 $ such that
\begin{equation}
\forall t\ge T, \quad \inf_{Q\in {\mathcal G}} \|u(t,\cdot)-\sum_{j=1}^N \lambda_j \,\exp
(\cdot-q_j) \|_{H^1} \le \gamma \label{ini4}
\end{equation}
and
\begin{equation}
\forall t\le -T, \quad \inf_{Q\in {\mathcal G}} \|u(t,\cdot)-\sum_{j=1}^N \lambda_{N+1-j}\, \exp
(\cdot-q_j) \|_{H^1} \le \gamma \; , \label{ini5}
\end{equation}
where $ {\mathcal G}:=\{Q\in \R^N, \, q_1<q_2<..<q_N\} $ and $ 0<\lambda_1<..<\lambda_N $ are the eigenvalues of the  matrix
 $ \Bigl( p_j^0 e^{-|q_i^0-q_j^0|/2}\Bigr)_{1\le i,j\le N}
$.
\end{cor}
This paper is organized as follows. In
Section~\ref{Well-posedness result} we state a well-posedness
result for (C-H) established in \cite{CM1} and \cite{dan1}. This
allows us to work in the function space $ Y([0,T]) $  that contains the peakons. Next, in Section~\ref{section1peak}
we present the result and the proof of Constantin and Strauss on
the stability of peakons. Section~\ref{sec-multi} is devoted to
the proof of Theorem~\ref{mult-peaks}. It is divided into  four
subsections. First we use a modulation argument  in order to
control the distance between the different bumps of the solution
we consider. Then we state a monotonicity result that was
established in~\cite{EM}. In Subsection~\ref{Localized energy
estimates} we establish a local version of an estimate involved in
the stability of a single peakon.
 The proof of  Theorem~\ref{mult-peaks} is completed in Subsection \ref{end-proof}.
In Section~\ref{sec-proof-cor} we recall some properties of
 the multipeakons
 and prove Corollary~\ref{cor-mult-peaks}. Finally in the appendix we give the proof of the monotonicity result for sake of completeness.

As mentioned above, the proof of the stability of trains of
peakons  does not enter the general framework (\cite{MMT},
\cite{EM}, \cite{EL}) on orbital stability of ordered trains of solitary
waves. However, the strategy of combining the orbital stability of
 a single solitary wave with a monotonicity  result seems to be quite
robust.

%%%%%%%%%%%%%%%%%%
%%%%%%%%%%%  Cauchy

\section{Well-posedness result}\label{Well-posedness result}
 Recall that the peakons do not belong to $ H^{3/2}(\R) $. To give a sens to these solutions,
 \re{CH} has to be rewritten as
 \begin{equation}
u_t -u_{txx}=- \frac{3}{2} \partial_x(u^2)-\frac{1}{2} \partial_x (u_x^2)
+\frac{1}{2}\partial_x^3 (u^2)  \label{CH2}
\end{equation}
or
 \begin{equation}
u_t +u u_x +(1-\partial_x^2)^{-1}\partial_x (u^2+u_x/2)=0  \label{CH3} \; .
\end{equation}
In \cite{CM1}, \cite{dan1} (see also \cite{L}) the following
existence and uniqueness result is
 derived.
\begin{theorem} \label{wellposedness}
Let $ u_0\in H^1(\R) $ with $ m_0:=u_0-u_{0,xx} \in {\cal M}(\R)
$ then there exists $ T=T(\|m_0\|_{\cal M})>0 $ and a unique
solution $ u \in Y([-T,T]) $ to (C-H)
with initial data $ u_0 $. The functionals $ E(\cdot) $
and $ F(\cdot) $ are constant along the trajectory and if $ m_0 $
has a definite sign then $ u $ is global in time. \\
Moreover, let $ \{u_{0,n}\}\subset H^1(\R) $ with $
\{u_{0,n}-\partial^2_x u_{0,n} \} $ bounded in $ {\cal M}_+ (\R) $
such that $ u_{0,n}\to u_0 $ in $ H^1(\R) $. Then, for all $ T> 0
$,
\begin{equation}
u_n \longrightarrow u \mbox{ in } C([-T,T]; H^1(\R)) \; .
\end{equation}
\end{theorem}
 Let
 us note that the last assertion of the above theorem is not explicitely contained in the works mentioned above. However, following the same arguments as those developped in these works (see for instance Section 5 of \cite{L}), one can prove that there exists a subsequence
  $ \{u_{n_k}\} $ of solutions of \re{CH} that converges in $ C([-T,T]; H^1(\R)) $ to some solution $ v$ of \re{CH} belonging to $ Y(-T,T) $. Since $ u_{0,n_k} $ converges to
   $ u_0 $ in $ H^1 $, it follows that $ v(0)=u_0 $ and thus $ v=u $ by uniqueness. This
    ensures that the whole sequence $ \{u_n\} $ converges to $ u $ in $C([-T,T]; H^1(\R))$
     and concludes the proof of the last assertion.

%%%%%%%%%%%%%%%%%%%
%%%%%%%%%%%% One peakon

\section{Stability of a single peakon}\label{section1peak}
Recall that the classical proof of orbital stability (see
\cite{Benjamin}, \cite{GSS}), successfully used in the case  $
\kappa>0 $ in \cite{CS2}, is based on the spectral properties of
the second differential operator of the invariant functional $ L_c(\cdot):=c
E(\cdot)-F(\cdot) $ evaluated at the solitary wave $ \varphi_c $.
Indeed, using a Liouville substitution, it can be shown that the
spectrum of the $ L^2 $-self-adjoint operator
$$
H_c:=L_c''(\varphi_{\kappa,c})=-\partial_x
\Bigl( (2c-2\varphi_{\kappa,c})\partial_x\Bigr)-6
\varphi_{\kappa,c}-2\partial_x^2  \varphi_{\kappa,c}+2(c-2\kappa)
$$
contains a unique negative eigenvalue which is simple and that $
0 $ is a simple eigenvalue associated with $ \partial_x
\varphi_{\kappa,c} $. The rest of the spectrum consists of a
finite number of positive eigenvalues and of the essential
spectrum  $ [2c-4\kappa, +\infty[ $. Therefore, controlling the
negative direction by modulating the velocity $ c $ and using that $ \langle
E'(\varphi_{\kappa,c}), u-\varphi_{\kappa,c}\rangle \sim 0 $
(since $ E(\cdot)$ is conserved) and the kernel direction by
choosing a suitable translation $ \varphi_{\kappa,c}(\cdot-r) $
of $ \varphi_{\kappa,c} $, the orbital stability is proven by
writing the Taylor expansion of $ c E(\cdot)-F(\cdot) $ at
$\varphi_{\kappa,c} $, recalling that $ c
E'(\varphi_{\kappa,c})-F'(\varphi_{\kappa,c}) $ vanishes.

Now, in the case $ \kappa=0 $, $ H_c $ is degenerate since $
\varphi_{\kappa,c}(0)=c $ and the Liouville substitution is no
more well-defined. However, Constantin and Strauss (cf.
\cite{CS1}) succeeded in proving the orbital stability by a direct
 approach (see also \cite{CM2} for another stability result using
Cazenave-Lions method). Actually, a by-product of their  proof is
the following very rigid property : for any function $ v $ in some $
H^1 $-neighborhood of $ \varphi_c $ it holds
$$
\|v-\varphi_c(\cdot-\xi)\|_{H^1}^2 \lesssim
|E(v)-E(\varphi_c)|+\sqrt{c\, |L_c(v)-L_c(\varphi_c)| } \; .
$$
where $ v(\xi) =\displaystyle \max_{\R} v $. Since $ E(\cdot) $
and $ F(\cdot) $ are conserved and are continuous functional on $
H^1(\R) $, this clearly leads to the orbital stability.

Since we will use similar considerations, we present here a sketch
of the proof of the stability of peakons (Theorem~\ref{1peakon})
proved by Constantin and Strauss in~\cite{CS1}.
\begin{theorem}\label{1peakon}
Let be given $c>0$ . There exist $C>0$ and $\varepsilon_0>0$ such
that  if $u\in C([0, T[ ; H^1(\R))$ is a solution of $\re{CH}$
 such that $ E(u(t)) $ and $ F(u(t)) $ are conserved quantities on $ [0,T[ $ and
   $\|u(0)-\varphi_c\|_{H^1}\leqslant\varepsilon^2$, then
\begin{equation}
\|u(t,\cdot)-\varphi_c(\cdot-r(t))\|_{H^1}\leqslant C
\sqrt{\varepsilon}, \quad \forall t\in [0,T[,
\end{equation}
where $r(t)\in\R$ is any point where the function $u(t,\cdot)$
attains its maximum.
\end{theorem}

\noindent The proof of this theorem is principally based on
the following lemma of $\cite{CS1}$.
\begin{lem}\label{1peakon-lemme}
For any $u\in H^1(\R)$ and $\xi\in\R$,
\begin{equation}\label{eq1}
E(u)-E(\varphi_c)= \|u-\varphi_c(\cdot-\xi)\|^2_{H^1}+4c(u(\xi)-c).
\end{equation}
For any $u\in H^1(\R)$, let $M=\max_{x\in\R}\{u(x)\}$, then
\begin{equation}\label{eq2}
F(u)\leqslant M E(u)-\frac{2}{3}M^3.
\end{equation}
\end{lem}
\begin{rem}
It is worth noticing that \re{eq1} ensures that the minimum of
the $ H^1 $-distance between $ u $ and $ \{ \varphi_c(\cdot-\xi),
\, \xi\in \R \} $ is exactly reached at any point $ \xi $ where $
u $ attains its maximum on $ \R $.
\end{rem}

 \noindent {\bf Proof of
Theorem\ref{1peakon}} Let $u\in C([0, T[ ; H^1(\R))$ be a solution
of $\re{CH}$ with
$\|u(0)-\varphi_c\|_{H^1}\leqslant\varepsilon^2$ and let
$\xi(t)\in\R$ be such that $u(t,\xi(t))=\max_{\R} u(t,\cdot) $.
 By the remark above, $ t\mapsto
 \|u(t)-\varphi_c(\cdot-\xi(t))\|_{H^1}$ is continuous on $ [0,T[ $
 and
 $\|u(0)-\varphi_c(\cdot-\xi(0))\|_{H^1}\leqslant\varepsilon^2$.
  Moreover, as shown in \cite{CS1}, it is no to hard to check that
   for any $v\in H^1(\R)$ such that
$\|u-\varphi_c\|_{H^1}<\gamma$ for some  $\gamma<1$, it holds
   \begin{equation}\label{eq3}
|E(u)-E(\varphi_c)|<4\, c  \gamma\textrm{ and }
|F(u)-F(\varphi_c)|<10\, c \gamma.
\end{equation}
>From the conservation laws it follows that for any $ t\in [0,T[ $
  \begin{equation}
|E(u(t))-E(\varphi_c)|<4\, c\varepsilon^2 \textrm{ and }
|F(u(t))-F(\varphi_c)|<10\, c \varepsilon^2 \; . \label{eq33}
  \end{equation}
 Therefore, by a classical continuity argument, it suffices to prove that for any
  $ v\in H^1(\R) $ satisfying \re{eq33} and
$ \|v-\varphi_c(\cdot-\xi)\|_{H^1}\leqslant\varepsilon^{1/4} $,
with $ v(\xi)=\max_{\R} v $,  it holds
 actually
 $$
 \|v-\varphi_c(\cdot-\xi)\|_{H^1}\lesssim \sqrt{\varepsilon} \quad
 .
 $$
 Setting $ M=v(\xi) $ and
 $\delta=c-M=c-v(\xi)$, we notice that  \re{eq1} ensures that for
  $\delta\leqslant 0$,
$$\|v-\varphi_c(\cdot-\xi)\|^2_{H^1}\leqslant E(u_0)-E(\varphi_c)
\lesssim \ep^2.$$ Hence to prove the stability it remains to
examine the case $\delta>0$, that is the maximum of the function
$u$ is less than the maximum of the peakon $\varphi_c$.
 Substituting $M$ by $c-\delta$ in $\re{eq2}$, using
 \re{eq33} and that
\begin{equation}\label{Ephi}
E(\varphi_c)=2 c^2 \textrm{ and } F(\varphi_c)=\frac{4}{3} c^3 \, ,
\end{equation}
  one can easily check that
$$
\frac{4}{3}c^3-O(\ep^2)\leqslant
(c-\delta)(2c^2+O(\ep^2))-\frac{2}{3}(c-\delta)^3 $$
which leads to
\begin{equation}
\delta^2(c-\delta/3)\leqslant O(\ep^2) \, . \label{tiko}
\end{equation}
On the other hand, on account of the hypothesis $ \|v-\varphi_c(\cdot-\xi)\|_{H^1}
\le \varepsilon^{1/4} $ and of the continuous embedding of $ H^1(\R) $ into $ L^\infty(\R) $, it holds $\delta<c/2$ for $ \varepsilon $ small enough.
 Therefore \re{tiko} ensures that
$\delta\leqslant C \ep$, the constant $C$ depending only on $c$.
This estimate on $\delta$ combining with  $\re{eq1}$ and
$\re{eq33}$ concludes the proof of Theorem~$\ref{1peakon}$.

%%%%%%%%%%%%%%%%%%%%%%%%

%%%%%%%%%%%% Proof of the result

\section{Stability of multipeakons}\label{sec-multi}
%Proof of Theorem~$\ref{mult-peaks}$}

For $ \alpha>0 $ and $ L>0 $ we define the following neighborhood of all the sums
  of N peakons of speed $ c_1,..,c_N $ with spatial shifts $ x_j $ that satisfied
$ x_j-x_{j-1}\ge L $.
  \begin{equation}
U(\alpha,L) = \Bigl\{
u\in H^1(\R), \, \inf_{x_j-x_{j-1}> L} \|u-\sum_{j=1}^N \varphi_{c_j} (\cdot-x_j) \|_{H^1} < \alpha \Bigr\}\; .
\end{equation}

 By the continuity of the map $ t\mapsto u(t) $ from $ [0,T[ $ into $ H^1(\R) $,
to prove Theorem \ref{mult-peaks} it suffices to prove that there exist
$ A>0 $, $ \varepsilon_0>0 $
and $ L_0>0 $ such  that $ \forall L>L_0 $ and $ 0<\varepsilon<\varepsilon_0 $, if
$u_0$ satisfies \re{ini} and if  for some $ 0<t_0< T $,
\begin{equation}\label{e11}
 u(t)\in U\left(A(\sqrt{\varepsilon}+L^{-{1/8}}),L/2\right)  \textrm{ for all }t\in[0,t_0]
\end{equation}
then
 \begin{equation}\label{e12}
 u(t_0) \in U\left(\frac{A}{2}(\sqrt{\varepsilon}+L^{-{1/8}}),\frac{2L}{3}\right).
\end{equation}
Therefore, in the sequel of this section we will assume \re{e11} for some
$ 0<\varepsilon<\varepsilon_0 $ and $ L>L_0 $, with $ A$, $ \varepsilon_0 $ and $ L_0 $ to be specified later, and we will prove
 \re{e12}.
%%%%%%%%%%%%%%%%%%%%%%%
%%%%%%%%%%%%% Modulation
\subsection{Control of the distance between the peakons}
In this subsection we want to prove that the different bumps of $ u $  that are individualy close to a peakon get away  from each others as time is increasing.
 This is crucial in our analysis since we do not know how to manage strong interactions.
\begin{lem}\label{eloignement}
Let $ u_0 $ satisfying \re{ini}. There exist $\alpha_0>0$,
$L_0>0$ and $C_0>0$  such that for all  $0<\alpha<\alpha_0$ and
$0<L_0<L$ if $u(t)\in U(\alpha, L/2) $ on $ [0,t_0] $ for some $
0<t_0< T $ then there exist $ C^1
$-functions  $ {\tilde x_1}, .., {\tilde x_N} $ defined on $
[0,t_0] $ such that
  \begin{equation}
  \frac{d}{d t} {\tilde x_i} = c_i +O(\sqrt{\alpha}) +O(L^{-1}) , \; i=1,..,N \, ,
  \label{vitesse}
  \end{equation}
  \begin{equation} \label{distH1}
  \|u(t)-\sum_{i=1}^N \varphi_{c_i} (\cdot -ž{\tilde x_i}(t)) \|_{H^1} =
  O(\sqrt{\alpha}) \, ,
  \end{equation}
\begin{equation}
{\tilde x_i}(t)-\tilde{x}_{i-1}(t) \ge 3L/4+(c_{i}-c_{i-1}) t/2 ,
 \; i=2,..,N . \label{eloi}
\end{equation}
Moreover, setting $ J_i:=[y_i(t), y_{i+1}(t)] $, $ i=1,..,N $,
 with
  \begin{equation} \label{defyi}
y_1=-\infty, \; y_{N+1}=+\infty \mbox{ and } y_i(t)=\frac{{\tilde x_{i-1}}(t)+{\tilde x_i}(t)}{2}\quad
  i=2,..,N,
  \end{equation}
  it holds
\begin{equation}
 |x_i(t)-{\tilde x_i}(t)| \le L/12 , \; i=1,..,N . \label{prox}
  \end{equation}
  where $ x_1(t),..,x_N(t) $ are any point such that
\begin{equation}
 u(t,x_i(t))=\max_{J_i(t)} u(t),  \; i=1,..,N .  \label{maxi}
  \end{equation}
\end{lem}
{\it Proof. }
To prove this lemma we use a modulation argument. The strategy is to
    construct
 $ N $ $ C^1 $-functions $ {\tilde x}_1,..,{\tilde x}_N $ on $ [0,t_0]$ satisfying
  a suitable orthogonality condition, see \re{mod2}.
   Thanks to this orthogonality condition we will be able to prove that the speed of the
    $ {\tilde x}_i $ stays close to $ c_i $ on $ [0,t_0] $.
\begin{rem}
It is crucial to note that in the previous works on stability of
sum of solitary waves (\cite{MMT}, \cite{EM}, \cite{EL}) one
needs similar modulation to ensure (among other things) that $v$
remains in a subspace of codimension two of $ H^1(\R) $ where the operator
 $ H_c $ (see the beginning of this section) is positive. Here, as already mentioned,
 we do not use such operator in the proof of orbital stability of
peakons but we still need a modulation to ensure that the
different bumps of $ u $ get away from each others.
\end{rem}
 For $Z=(z_1,..,z_N)\in \R^N $ fixed such that $ z_i-z_{i-1}>L/2 $, we set
$$
R_Z(\cdot)=\sum_{i=1}^N \varphi_{c_i}(\cdot -z_i) \quad.
$$
 For  $0<\alpha<\alpha_0 $ we define the function
 \begin{eqnarray*}
 Y\, :\, (-\alpha,\alpha)^N \times B_{H^1}(R_Z,\alpha)
 & \to &  \R^n \\
 (y_1,..,y_N,u) & \mapsto & (Y^1(y_1,..,y_N,u), ..,Y^N(y_1,..,y_N,u))
 \end{eqnarray*}
 with
 $$
  Y^i(y_1,..,y_N,u)
=\int_{\R} \Bigl(u-\sum_{j=1}^N \varphi_{c_j}(\cdot-z_j-y_j)\Bigr)\partial_x
\varphi_{c_i}(\cdot-z_i-y_i) \; .
 $$
 $Y $ is clearly of class $ C^1 $. For $ i=1,..,N $,
\begin{equation}
\frac{\partial Y^i}{\partial y_i}(y_1,..,y_N,u)  =
\int_{\R}  \Bigl(u_x- \sum_{j=1, j\neq i }^N
\int_{\R}\partial_x \varphi_{c_j}(\cdot-z_j-y_j)\Bigr) \partial_x
\varphi_{c_i}(\cdot-z_i-y_i)
 \, dx  \, .
\end{equation}
 and  $\forall j\neq i $
 $$
 \frac{\partial Y^i}{\partial y_j}(y_1,..,y_N,u) = \int_{\R}
 \partial_x \varphi_{c_j}(\cdot-z_j-y_j) \partial_x
\varphi_{c_i}(\cdot-z_i-y_i)
 \, dx \; .
 $$
Hence,
\begin{eqnarray}
\frac{\partial Y^i}{\partial y_i}(0,..,0,R_Z) &  =  &
\|\partial_x \varphi_{c_i}\|_{L^2}^2\ge c_1^2 \, .
\end{eqnarray}
and, for $ j\neq i $, using the exponential decay of $ \varphi_c $  and that $ z_i-z_{i-1}>L $ we infer that
 for $ L_0 $ large enough (recall that $L>L_0$),
\begin{eqnarray*}
\frac{\partial Y^i}{\partial y_j}(0,..,0,R_Z) &  =  &
\int_{_R} \partial_x \varphi_{c_j}(\cdot-z_j)\,  \partial_x
 \varphi_{c_i}(\cdot-z_i) \, dx
 \nonumber \\
& \le &  O(e^{-L/4}) \;.
\end{eqnarray*}
We deduce that, for $ L >0 $ large enough, $ D_{(y_1,..,y_N)} Y(0,..,0,R_Z) =D+P $ where $ D $ is an invertible diagonal matrix with $ \|D^{-1}\|
\le (c_1)^{-2} $ and $  \|P\|\le O(e^{- L/4}) $. Hence
there exists $L_0> 0 $ such that for $ L>L_0 $, $ D_{(y_1,..,y_N)} Y(0,..,0,R_Z) $ is invertible with an inverse matrix
 of norm smaller than 2 $(c_1)^{-2} $. From the implicit function theorem we deduce that there exists $ \beta_0>0 $
  and $ C^1 $ functions $(y_1,..,y_N) $ from $ B(R_Z,\beta_0) $ to  a neighborhood
  of $ (0,..,0)$ which are uniquely determined such that
  $$
  Y(y_1,..,y_N,u) = 0 \mbox{ for all } u \in B(R_Z, \beta_0) \; .
  $$
  In particular, there exits $ C_0>0 $ such that if $ u\in  B(R_Z,\beta) $,
   with $ 0<\beta \le \beta_0 $, then
  \begin{equation}
  \sum_{i=1}^N |y_i(u)|\le C_0 \beta \ ;. \label{estyi}
  \end{equation}
  Note that $ \beta_0 $ and $ C_0 $ only depend on $ c_1 $ and $ L_0 $ and not on the
   point $ (z_1,..,z_N) $.
   For $ u\in  B(R_Z,\beta_0) $
   we set $ {\tilde x}_i(u)=z_i+y_i(u) $. Assuming that $ \beta_0\le \frac{L_0}{8C_0} $,
   $({\tilde x}_1,..,{\tilde x}_N)$ are thus $ C^1 $-functions on $  B(R_Z,\beta) $
satisfying
   \begin{equation}
{\tilde x}_j(u)-{\tilde x}_{j-1}(u) >  L/2-2C_0\beta> L/4 \quad
. \label{opop}
\end{equation}
For $ L\ge L_0 $ and $ 0<\alpha<\alpha_0< \beta_0/2 $ to be chosen later, we
define the modulation of $ u \in U(\alpha, L/2 )$
 in the following way : we cover the trajectory of $ u$
 by a finite number of open balls in the following way :
$$
\Bigl\{u(t),\, t\in [0,t_0]\Bigr\} \subset
 \displaystyle\bigcup_{k=1,.., M }
  B(R_{Z^k},2\alpha)
 $$
 It is worth noticing
that, since $ 0<\alpha<\alpha_0< \beta_0/2 $,
 the functions  $ {\tilde x}_j(u) $ are uniquely
  determined for $ u \in B(R_{Z^k)},2 \alpha) \cap B(R_{Z^{k'}},2\alpha) $.
We can thus define the functions $ t \mapsto {\tilde x}_j(t) $ on
  $ [0,t_0] $ by setting
 $ {\tilde x}_j(t)={\tilde x}_j(u(t)) $.
By construction
\begin{equation}
\int_{\R} \Bigl( u(t,\cdot) -\sum_{j=1}^N \varphi_{c_j}(\cdot-{\tilde x_j}(t)) \Bigr)
 \partial_x \varphi_{c_i} (\cdot -{\tilde x_i}(t)) \, dx = 0 \; . \label{mod2}
\end{equation}
Moreover, on account of  \re{estyi} and the fact that $
\varphi_c''$ is the sum of a $ L^1 $ function and a Dirac mass  it
holds
\begin{equation}
\|v(t)\|_{H^1} \lesssim C_0 \sqrt{\alpha}\, , \quad \forall t\in
[0,t_0] \; . \label{estepsilon}
\end{equation}
Let us now prove that the speed of $ {\tilde x}_i $ stays close to $ c_i $.
 We set
 $$
 R_j(t)=\varphi_{c_j}(\cdot-{ \tilde x}_j(t)) \mbox{ and }
 v(t)=u(t)-\sum_{i=1}^N R_j(t)=u(t,\cdot)-R_{{\tilde X}(t)} \; .
$$
Differentiating \re{mod2} with respect to time we get
$$
\int_{\R} v_t  \partial_x R_i =\dot{\tilde x}_i \, \langle
 \partial_x^2 R_i \, ,\,  v  \rangle_{H^{-1}, H^1}, \; .
$$
and thus
\begin{equation}
\Bigl|\int_{\R} v_t  \partial_x R_i\Bigr|
\le |\dot{\tilde x}_i| O(\|v\|_{H^1}) \le |\dot{\tilde x}_i-c_i|
 O(\|v\|_{H^1})+O(\|v\|_{H^1})\; . \label{huhu}
\end{equation}
Substituting $ u $ by $ v+\sum_{j=1}^N R_j $ in \re{CH3}  and
using that $ R_j $
 satisfies
 $$
 \partial_t R_j +(\dot{\tilde x}_j-c_j)  \partial_x R_j + R_j \partial_x R_j
 +(1-\partial_x^2)^{-1}\partial_x  [u^2+u_x^2/2] = 0 \;,
 $$
 we infer that $ v$ satisfies on $ [0,t_0] $,
  \arraycolsep1pt
 \begin{eqnarray}
 v_t&  - & \sum_{j=1}^N (\dot{x}_j-c_j)  \partial_x R_j
= -\frac{1}{2}  \partial_x \Bigl[(v+\sum_{j=1}^N
 R_j)^2- \sum_{j=1}^N R_j^2  \Bigr] \nonumber \\
 & &-(1-\partial_x^2)^{-1}\partial_x \Bigl[(v+\sum_{j=1}^N R_j)^2- \sum_{j=1}^N R_j^2
 +\frac{1}{2} (v_x +\sum_{j=1}^N \partial_x R_j)^2 -\frac{1}{2}\sum_{j=1}^N
  (\partial_x R_j)^2\Bigr]\; . \nonumber
 \end{eqnarray}
 \arraycolsep5pt
Taking the $ L^2 $-scalar product with $ \partial_x R_i $, integrating by parts, using the
 decay of $ R_j $ and its first derivative, \re{estepsilon}, \re{huhu} and \re{opop}, we find
 \begin{equation}
 |\dot{\tilde x}_i-c_i|\Bigl(\|\partial_x R_i \|_{L^2}^2 +O(\sqrt{\alpha}) \Bigr)
 \le O(\sqrt{\alpha}) +O(e^{ -L/8})\; . \label{fofo}
 \end{equation}
Taking $\alpha_0$ small enough and $ L_0 $ large enough we
  get
$ |\dot{\tilde x}_i-c_i| \le (c_i-c_{i-1})/4 $ and thus  for all $ 0<\alpha<\alpha_0
$ and $ L\ge L_0>3C_0\varepsilon $, it follows from \re{ini}, \re{estyi} and \re{fofo}   that
\begin{equation}
{\tilde x}_j(t)-{\tilde x}_{j-1}(t)>
L-C_0\varepsilon+(c_j-c_{j-1}) t/2 , \quad \forall t\in[0,t_0]
\; . \label{xj-xj-1}
\end{equation}
which yields \re{eloi}.\\
Finally from \re{estepsilon} and the continuous embedding of $
H^1(\R)$ into $ L^\infty(\R) $, we infer that
$$ u(x) = R_{\tilde
X}(x)+O(\sqrt{\alpha}),  \quad \forall x\in \R \, .
$$
Applying this formula with $ x=x_i=\max_{J_i(t)} u(t) $ and taking advantage of
 \re{eloi}, we obtain
$$
u(x_i)
=c_i+O(\sqrt{\alpha})+O(e^{-L/4}) \ge
2c_i/3 \; . $$
 On the other hand, for $ x\in J_i\backslash
]{\tilde x_i}-L/12,{\tilde x_i}+L/12[ $, we get
$$ u(x)\le c_i
e^{-L/12}+O(\sqrt{\alpha})+O(e^{-L/4}) \le c_i/2 \; . $$ This
ensures that  $x_i$ belongs to $ [{\tilde x_i}-L/12,{\tilde
x_i}+L/12] $.
%%%%%%%%%%%%%%%%%%%%%%%%

%
%%%%%%%%%%%%%%%%%%%%%%%%
%%%%%%%%%%%% Monotonie

\subsection{Monotonicity property}\label{Sectmonotonie}
Thanks to the preceding lemma, for $ \varepsilon_0> 0 $ small enough and $ L_0>0 $ large enough, one can construct  $C^1$-functions $ {\tilde x_1}, .., {\tilde x_N} $ defined on $
[0,t_0] $ such that \re{vitesse}-\re{prox} are satisfied.
In this subsection we state the almost monotonicity of  functionals that are very close
 to the energy at the right  of the $ i $th bump, $ i=1,..,N-1 $ of $ u $. The proof
   is  similar to the one of
  Lemma~4.2 in~\cite{EL}. We give it in the appendix for sake of completeness.\\
%However, since we do not try this time to get an optimal decay rate as
%   in Theorem 3.1,
 Let $ \Psi $ be a  $ C^\infty $ function such that
 $ 0<\Psi\le 1 $,  $ \Psi'>0 $ on $ \R $, $ |\Psi'''|\le 10 |\Psi'| \mbox{ on } [-1/2,1/2]
 $,
 $$
 \Psi(x)=\left\{ \begin{array}{ll}
 e^{-|x|} & \quad x<-1/2 \\
 1-e^{-|x|}& \quad x>1/2
 \end{array}
 \right.\quad \mbox{ and }
\left\{ \begin{array}{ll}
 \Psi(x)\le 2e^{-|x|}& \mbox{ on } [-1/2,0]\\
 1-\Psi(x) \le 2e^{-|x|} & \mbox{ on } [0,1/2]
\end{array}
\right.
$$
 Setting $ \Psi_K=\Psi(\cdot/K) $, we introduce for $ j\in \{2,..,N\}, $
 $$
 I_{j,K}(t)=I_{j,K}(t,u(t))=  \int_{\R} (u^2(t)+u_x^2(t)) \Psi_{j,K}(t) \, dx\,,
 $$
 where $ \Psi_{j,K}(t,x)=\Psi_K(x-y_j(t)) $ with $ y_j(t)$,
 $ j=2,..,N $, defined in \re{defyi}. Note that $ I_j(t) $ is close to $  \|u(t)\|_{H^1(x>y_j(t))} $ and thus measures the energy at the right of the
   $ (j-1)$th bump of $ u$.
   Finally, we set
   \begin{equation}
\sigma_0=\frac{1}{4} \min \Bigl(c_1,c_2-c_1,..,c_N-c_{N-1}\Bigr)
\; .
   \end{equation}
   In \cite{EL} the following monotonicity result is derived.
\begin{lem}\label{monotonicitylem}
  Let $ u\in Y([0,T[) $ be a solution of (C-H) satisfying \re{distH1} on
   $[0,t_0] $.
 There exist $ \alpha_0>0 $ and $ L_0>0 $ only depending on $ c_1 $ such that
    if $ 0<\alpha<\alpha_0 $ and $ L\ge L_0 $ then for any $ 4\le K \lesssim L^{1/2} $,
    \begin{equation}\label{monotonicityestim}
    I_{j,K}(t)-I_{j,K}(0)\le O( e^{-\frac{\sigma_0 L }{8K}}) , \quad \forall j\in\{2,..,N\}, \; \quad \forall t\in [0,t_0] \; .
    \end{equation}
\end{lem}

%%%%%%%%%%%%%%%%%%%
%%%%%%%%%%%% energie

\subsection{A localized and a global estimate}\label{Localized energy estimates}

We define the function $ \Phi_i=\Phi_i(t,x) $  by $
\Phi_1=1-\Psi_{2,K}=1-\Psi_K(\cdot-y_2(t))$, $
\Phi_N=\Psi_{N,K}=\Psi_K(\cdot-y_N(t)) $ and for $i=2,..,N-1 $
$$
\Phi_i=\Psi_{i,K}-\Psi_{i+1,K}=\Psi_K(\cdot-y_i(t))-\Psi_K(\cdot-y_{i+1}(t))\; ,
$$
where $ \Psi_K $ and the $ y_i $'s are defined in Section 4.2.
It is easy to check that $\displaystyle \sum_{i=1}^N \Phi_{i,K}\equiv 1 $.
 We take $L>0$ and $L/K>0$ large enough so that  $ \Phi_i $ satisfies
 \begin{equation}
 |1-\Phi_{i,K}| \le  4 e^{-\frac{L}{4K}} \mbox{ on } [\widetilde x_i-L/4, \widetilde x_i+L/4]
  \label{de1}
 \end{equation}
 and
\begin{equation}
 |\Phi_{i,K}| \le 4 e^{-\frac{L}{4K}} \mbox{ on } [\widetilde x_j-L/4,
 \widetilde x_j+L/4] \mbox{ whenever } j\neq i \; .\label{de2}
 \end{equation}
 We will use  the following localized  version of $ E $ and $ F $ defined  for 
$i\in \{1,..,N\}, $ by
 \begin{equation}
 E_i^t(u) = \int_{\R} \Phi_{i}(t) (u^2+u_x^2) \mbox{ and }
 F_i^t(u)= \int_{\R} \Phi_i(t) (u^3+u u_x^2) \; .
 \end{equation}
{\bf Please note that henceforth we take $K=L^{1/2}/8 $.} \\
The following lemma gives a localized version of \re{eq2}.
Note that the functionals $ E_i $ and $ F_i $ do not depend on
time in the statement below  since we fix $ {\tilde
x_1}<..<{\tilde x_N} $.
\begin{lem}\label{m-p-lemme}
 Let  be given $ N $ real numbers $ {\tilde x_1}<..<{\tilde x_N} $
  with $ {\tilde x_i}-{\tilde x_{i-1}} \ge 2L/3 $.
 Define the $ J_i $'s as in \re{defyi} and assume that, for $ i=1,..,N $, there exists
  $ x_i\in J_i $ such that $ |x_i-{\tilde x_i}|\le L/12 $ and $ u(x_i)=\displaystyle \max_{J_i} u:=M_i $.
  Then, for any $ u\in H^1(\R) $,  it holds
\begin{equation}\label{eq2m}
F_i(u)\leqslant M_i E_i(u)-\frac{2}{3}M_i^3+\|u_0\|_{H^1}^3 O(L^{-{1/2}}), \quad
 i\in \{1,..,N\} \, .
\end{equation}
\end{lem}
{\it Proof. } Let $ i\in \{1,..,N\}$ be fixed. Following \cite{CS1},
we introduce the function  $g$ defined  by
$$g(x)=\left\{
\begin{array}{l}
u(x)-u_x(x) \; \mbox{ for } \; x<{ x_i} \\
u(x)+u_x(x) \; \mbox{ for }\; x>{ x_i}
\end{array}.
\right.
$$
Integrating by parts we compute
\begin{eqnarray}\label{ug2}
\int ug^2\Phi_i&=&\int_{-\infty}^{ x_i}(u^3+uu_x^2-2u^2u_x)\Phi_i
+\int_{ x_i}^{+\infty}(u^3+uu_x^2+2u^2u_x)\Phi_i\nonumber\\
&=&F_i(u)-\frac{4}{3}u({x_i})^3\Phi_i({
x_i})+\frac{2}{3}\int_{-\infty}^{ x_i}u^3\Phi_i^{'}
-\frac{2}{3}\int_{ x_i}^{+\infty}u^3\Phi_i^{'}\, .
\end{eqnarray}
Recall that we take $ K=\sqrt{L}/8 $ and thus  $ |\Phi'|\le C/K = O(L^{-1/2}) $.
 Moreover, since $ |x_i-{\tilde x_i}|\le L/12 $, it follows from \re{de1} that     $ \Phi-i({ x
_i})=1+O(e^{-L^{1/2}}) $  and thus
\begin{equation}
\int u g^2 \Phi_i = F_i(u)-\frac{4}{3} M_i^3+\|u_0\|_{H^1}^3 O(L^{-1/2}) \; .
\end{equation}

On the other hand,
\begin{eqnarray}\label{hg2}
\int u g^2\Phi_i& \le & M_i \int g^2 \Phi_i \nonumber \\
&\leq &M_i \Bigl( E_i(u)-2\int_{-\infty}^{x_i}uu_x\Phi_i+2\int_{
x_i}^{+\infty}
uu_x\Phi_i\nonumber\\
&\leq &M_iE_i(u)-2M_i^3+\|u_0\|_{H^1}^3 O(L^{-1/2}) \; .
\end{eqnarray}
This proves \re{eq2m}.\vspace{2mm} \\
Now let us state a global identity related to \re{eq1}.
\begin{lem}
For any $ Z\in \R^N $ such that $ |z_i-z_{i-1}|\ge L/2 $ and any $ u\in H^1 $ it holds
\begin{equation} \label{rf}
 E(u)-\sum_{i=1}^N E(\varphi_{c_i}) =
 \| u-R_Z \|_{H^1}^2 +4\sum_{i=1}^N c_i (u(z_i)-c_i)+ O(e^{-L/4}) \; .
\end{equation}
\end{lem}
{\it Proof .} Using the relation between $ \varphi $ an its derivative and  integrating  by parts, we get
\begin{eqnarray*}
E(u-R_Z) & = & E(u)+E(R_Z) -2\sum_{i=1}^N \int u\,  \varphi_{c_i}(\cdot-z_i) + u_x \, \partial_x \varphi_{c_i}(\cdot-z_i)   \\
& =& E(u)+E(R_Z) -2\sum_{i=1}^N \int u \, \varphi_{c_i}(\cdot-z_i) \\
 & & +2 \sum_{i=1}^N
 \int_{z_i}^{+\infty} u_x \, \varphi_{c_i}(\cdot-z_i) -2\sum_{i=1}^N
 \int^{z_i}_{-\infty} u_x \, \varphi_{c_i}(\cdot-z_i)  \\
 & = & E(u)+E(R_Z) -4 \sum_{i=1}^N c_i u(z_i) \; .
 \end{eqnarray*}
On the other hand, since $ |z_i-z_{i-1}|\ge L/2 $, it is not too hard to check that
$$
E(R_Z)= \sum_{i=1}^N E(\varphi_{c_i}) +O(e^{-L/4}) = 2\sum_{i=1}^N c_i^2+O(e^{-L/4})\, .
$$
Combining these two identity, the desired result follows. \vspace*{2mm} \\
As a consequence of this lemma, we obtain an estimate on the $ H^1 $ distance
 between $ u(t) $ and $ R_{X(t)} $.
\begin{lem}
Under the same hypotheses as in Lemma \ref{eloignement}, the function $ X=(x_1,..,x_N) $ constructed in Lemma \ref{eloignement} satisfies on $ [0,t_0 ]$,
 \begin{equation}
\|u(t)-R_{X(t)} \|_{H^1} \le O(\alpha)+O(e^{-L/8}) \; .\label{ri}
 \end{equation}
\end{lem}
{\it Proof. }
Since $ u(t)\in U(\alpha, L/2) $ for $ t\in [0,t_0] $, on account of Lemma \ref{eloignement} for any $ t\in [0,t_0] $ there exists $ Z=(z_1,..,z_n) $ with $
 z_i\in J_i(t) $ such that $ E(u(t)-R_Z) = O(\alpha^2) $. Recalling that  $ u(t,x_i(t))=\max_{J_i(t)} u(t) $, we deduce \re{ri}  from \re{rf}.

%%%%%%%%%%%%%%%%%%%%%%%%%%%%%%%%%
\subsection{End of the proof of Theorem~$\ref{mult-peaks}$}\label{end-proof}

Recall that $\sum_{i=1}^{N}E_i(v)=E(v)$ for any $v\in H^1(\R)$. From \re{ini}
 it is easy to check that
\begin{equation}
E(u(t))=E(u_0)=\sum_{j=1}^N E(\varphi_{c_j})+O(\varepsilon^2) +O(e^{-L/4}) , \;
 \forall t\in [0,T] \; .\label{fu}
\end{equation}
Let us set $M_i=u(t_0,x_i(t_0)) $ and $ \delta_i=c_i-M_i$. To
conclude the proof, it thus suffices to prove that there exists $C>0$
which does not depend on $A$ such that
\begin{equation}\label{deltai}
\delta_i\leq C (\ep+L^{-{1/4}}) \textrm{ for all } i.
\end{equation}
Indeed, in this case \re{fu} and \re{rf}, with $ Z=X(t_0)$, ensure the existence of $C>0$
independent of $A$ such that
$$\|u-\sum_{j=1}^N \varphi_{c_j} (\cdot- x_j) \|_{H^1} < C (\ep^{1/2}+L^{-{1/8}}),$$
so that one can take $A=2C$ to conclude the proof (Recall that we already know
 from \re{eloi}-\re{prox} that $ x_i-x_{i-1} \ge 2L/3 $ for $ i\in\{2,..,N\} $).
Let us prove \re{deltai}. From~\re{eq2m} by taking the sum over $i$ one gets :
$$F(u(t_0))=\sum_{i=1}^{N}F_i(u(t_0))\leqslant \sum_{i=1}^{N}M_i E_i(u(t_0))-\frac{2}{3}\sum_{i=1}^{N}M_i^3+O(L^{-{1/2}})$$
Setting  $ \Delta_0^{t_0} F(u)=F(u(t_0))-F(u(0)) $ and
 $ \Delta_0^{t_0} E(u)=E(u(t_0))-E(u(0)) $, this implies
\begin{eqnarray}\label{e1}
0=\Delta_0^{t_0} F(u)=\sum_{i=1}^{N}\Delta_0^{t_0} F_i(u)&\leqslant& \sum_{i=1}^{N}M_i\Delta_0^t E_i(u)
-2/3 \sum_{i=1}^{N} M_i^3 \\
&&+\sum_{i=1}^{N}(-F_i(u_0)+M_i E_i(u_0))+O(L^{-{1/2}})\nonumber
\end{eqnarray}
By \re{ini}, the exponential decay of the $ \varphi_{c_i} $'s   and the $ \Phi_i $'s, and the definition
 of $ E_i $ and $ F_i $, it is easy to check that
 $$
|E_i(u_0)-E(\varphi_{c_i}|+
|F_i(u_0)-F (\varphi_{c_i}|\le O(\varepsilon^2)+O(e^{-\sqrt{L}}), \; \forall i\in\{1,..,N\}\, .
$$
Setting $M_0=0$ and using \re{Ephi}, one thus finds after
 having substituted  $M_i$ by $c_i-\delta_i $ that
\begin{equation}\label{e2}
\sum_{i=1}^{N}(-F_i(u_0)+M_i E_i(u_0)-2/3 M_i^3)
=2\sum_{i=1}^{N}(-c_i\delta_i^2+\frac{1}{3}\delta_i^3)+O(\ep^2)
+O(e^{-\sqrt{L}})\,
.
\end{equation}
Note that by \re{ri} and the continuous embedding of $ H^1(\R) $ into $ L^\infty(\R) $,
 $ M_i=c_i+O(\alpha)+O(e^{-L/8})$, and thus
 \begin{equation}
0<M_1<\cdot \cdot<M_N \mbox{ and } \delta_i<c_i/2  \label{zq}
 \end{equation}
 for $ \alpha_0=A(\sqrt{\varepsilon_0}+L_0^{-1/8}) $ small enough.
 Using   the Abel transformation and the monotonicity
estimates~\re{monotonicityestim}, we thus get
\begin{eqnarray}\label{e3}
\sum_{i=1}^{N}M_i\Delta_0^t
E_i(u)&=&\sum_{i=1}^{N}(M_i-M_{i-1})\Delta_0^t I_i
\leqslant O(\ep^2+e^{- \sqrt{L}})\; .
\end{eqnarray}
Injecting $(\ref{e2})$ and $(\ref{e3})$ in $(\ref{e1})$ we obtain
\begin{equation}\label{e4}
\sum_{i=1}^{N}(c_i\delta_i^2-\frac{1}{3}\delta_i(t)^3)
=\sum_{i=1}^{N}\delta_i^2(c_i-\frac{1}{3}\delta_i)\leqslant O(\ep^2+L^{-{1/2}}).
\end{equation}
\re{zq} and \re{e4} yield \re{deltai}
and concludes the proof of the theorem.

\section{Proof of Corollary~\ref{cor-mult-peaks}}\label{sec-proof-cor}
As written in the introduction, Camassa and Holm discovered that \re{CH} possesses special solutions given by
\begin{equation}
u(t,x)=\sum_{i=1}^N p_i(t) e^{|x-q_i(t)|} \label{multipeakons2}
\end{equation}
where the $ (p_i,q_i)\in (\R^2) $ satisfy the Hamiltonian system
\begin{equation} \label{systH}
\Bigl\{
\begin{array}{l}
\dot{q}_i=\sum_{j=1}^N p_j e^{-|q_i-q_j|} \\
\dot{p}_i = \sum_{j=1}^N p_i p_j \sgn (q_i-q_j) e^{-|q_i-q_j|}\; .
\end{array}
\Bigr.
\end{equation}
It is easy to check that the local solution of this differential system can be extended
 as soon as the $ q_i's $ stay distinct from each other.
In \cite{HR}, Holden and Raynaud proved that this is indeed the case if at time $ t=0 $, the $ p_i $ are all positive , i.e. there are only peakons (the case with only anti-peakons works also but in the case with peakon and anti-peakon this is no longer true).
 More precisely,
 they  proved that if at time $ t= 0 $,
\begin{equation}
\label{initH}
p_1,.., p_N >0 \mbox{ and } q_1<q_2<q_N
\end{equation}
then \re{initH} remains true for all  time. In particular, under these hypotheses the different peakons never overlap each others. For example, if a larger peakon
 follows a smaller one, it will come close to this last one and then transfer part
  of its energy to it. In this way, the smaller one will become the larger one
   and the two peakons will be well ordered.
In \cite{Beals} (see also \cite{CH1}), using the integrability of \re{CH}, Beals {\it et al} established a formula for the asymptotics of the $ q_i$'s and the $ p_i$'s.
 In particular, they prove the following limits for the $ p_i $ and
  $ {\dot q_i}$, $i\in\{1,..,N\}$,
 \begin{equation}
\lim_{t\to +\infty}  p_i(t)=\lim_{t\to +\infty} {\dot q_i}(t) =\lambda_i \label{limit1}
 \end{equation}
and
\begin{equation}
\lim_{t\to-\infty}  p_i(t)=\lim_{t\to -\infty} {\dot q_i}(t) =\lambda_{N+1-i}\; ,\label{limit2}
 \end{equation}
where $ 0<\lambda_1<\cdot \cdot <\lambda_N $ are the eigenvalues of the matrix
 $(p_j(0) e^{-|q_i(0)-q_j(0)|/2})_{i,j} $.
 \begin{rem}
 The matrix $A_N:=(p_j e^{-|q_i-q_j|/2})_{1\le i,j\le N} $ is obtained by substituing the multipeakon solution \re{multipeakons2}
  in the isospectral problem
  \begin{equation}
  \Psi_{xx}=\Bigl(\frac{1}{4}-\frac{m(t,\cdot)}{2\lambda} \Bigr) \Psi , \quad
   \mbox{ with } \, m=u-u_{xx}, \label{isospectral}
  \end{equation}
  associated with the Camassa-Holm equation. More precisely, any solution of \re{isospectral} with
   $ m=2\sum_{i=1}^N p_i \delta_{q_i}$, that vanishes at $ \mp \infty $, is completely determined by
    its values at the $q_j $'s and satisfies
   \begin{equation}
    \lambda \Psi(q_i)=\sum_{j=1}^N p_j e^{-|q_i-q_j|/2} \Psi(q_j), \quad \forall i\in \{1,..,N\}\,.
    \label{iso2}
    \end{equation}
In \cite{Beals}, \re{isospectral} is transformed into a density problem on $ [-1,1] $ by applying a Liouville transformation. The corresponding N-multipeakon matrix is then proved to possess  $ N $
 distinct positive eingenvalues. The arguments of \cite{Beals} hold also clearly for $ A_N $. Indeed, first since for any fixed $ \lambda $, \re{isospectral} has clearly at most one
   solution (up to multiplication by a scalar) that vanishes at $ \mp \infty $, it follows that the
  eigenvalues of $ A_N $ are all of geometric multiplicity one. Next, setting $ D=\mbox{diag}\, (p_i)$
   and $\Lambda_{i,j}=e^{-|q_i-q_j|/2} $,  $ A_N $ can be rewritten as $ D \Lambda $.
   Since $ \Lambda $ is symmetric with $ \Lambda_{ii}=1 $ and $ |\Lambda_{ij}|<1 $ for $ i\neq j $,
    $  \Lambda $ is actually  positively defined. Therefore there exists $ B $ a symmetric positively defined matrix such that $ \Lambda=B^2 $. It is then easy to check that $ A_N $ and $ BDB $ have got the
     same spectrum  and since
     $ BDB $ is symmetric  positively defined, this ensures that $ A_N$ possesses $ N $ distinct positive eigenvalues.
 \end{rem}

  Now, let be given  $ (p_i(0), q_i(0)) $ satisfying \re{initH} and $ \gamma>0 $. From the asymptotics above
   there exists $ T>0 $ such that
   \begin{equation}
   \label{fy}
 q_i(T)-q_{i-1}(T)>L       \mbox{ and } q_{i}(-T)-q_{i-1}(-T)>L
  \end{equation}
  with
  \begin{equation}
L>\max\Bigl( L_0,(\frac{\gamma}{2A})^8\Bigr) \; .
  \end{equation}
   From the last assertion
   of Theorem \ref{wellposedness}, for any given $ B>0 $, there exists $ \alpha> 0 $ such that if $ u_0 $ satisfies \re{ini3} then for all $ t\in[-T,T] $,
  \begin{equation} \label{fy2}
\Bigl\|u(t)-\sum_{i=1}^N p_i(t) e^{|x-q_i(t)|}\Bigr\|_{H^1} \le \Bigl(
\frac{\gamma}{2A}\Bigr)^4 \; .
   \end{equation}
   At this stage, it is crucial to remark that since \re{CH} is invariant under
    the transformation $ (t,x) \mapsto (-t,-x) $, Theorem \ref{mult-peaks} remains
     true when replacing $ t$ by $ -t $, $ z_j^0 $ by $ -z_j^0 $ and $ x_j(t) $
      by $ -x_j(-t) $. This gives a stability result in the past for trains of
       peakons that are ordered in the inverse order with respect to
        Theorem \ref{mult-peaks}.\\
Combining \re{fy}, \re{fy2}, Theorem \ref{mult-peaks} and the remark above, the first part of the corollay follows.

Finally, from \re{limit1}-\re{limit2}, we can also assume that
$$
|p_i(T)-\lambda_i|\le \frac{1}{100N} \Bigl(\frac{\gamma}{2A}\Bigr)^4 \mbox{ and }
|p_i(-T)-\lambda_{N-i}|\le \frac{1}{100N} \Bigl(\frac{\gamma}{2A}\Bigr)^4
$$
so that
$$
\Bigl\|u(T)-\sum_{i=1}^N \lambda_i e^{|x-q_i(T)|}\Bigr\|_{H^1} \le \Bigl(
\frac{\gamma}{2A}\Bigr)^4 \mbox{ and }
 \Bigl\|u(-T)-\sum_{i=1}^N \lambda_{N-i}\,  e^{|x-q_i(-T)|}\Bigr\|_{H^1} \le \Bigl(
\frac{\gamma}{2A}\Bigr)^4\; .
$$
This completes the proof of the corollary.

\section{Appendix}
{\it Proof of Lemma \ref{monotonicitylem}. }
Let us assume that $ u $ is smooth since the case $ u\in  Y([0,T[) $ follows by modifying slightly the arguments (see Remark 3.2 of \cite{EM}). From \re{CH2}, it is not too hard
 to check that for any smooth space function $ g  $, the  folllowing differential
  identity on the weighted energy holds :
\begin{eqnarray}
\frac{d}{dt}\int_{\R} (u^2+u_x^2) g \, dx&=&\int_{\R}(u^3+4uu_x^2)g^{'} \, dx\nonumber\\
&&\hspace*{-10mm}
-\int_{\R}u^3g^{'''} \, dx-\int_{\R}ug^{'}(1-\partial_x^2)^{-1}(2u^2+u_x^2) \, dx.
\label{go}
\end{eqnarray}
Applying \re{go} with $ g=\Psi_{j,K} $ one gets
\begin{eqnarray}
\frac{d}{dt}\int_{\R}  \Psi_{j,K}(u^2+u_x^2)  \, dx&=&-\dot{y_j}
\int_{\R} \Psi_{j,K}' (u^2+u_x^2) +\int_{\R}\Psi_{j,K}' (u^3+4uu_x^2) \, dx\nonumber\\
&&-\int_{\R}\Psi_{j,K}^{'''} u^3\, dx-\int_{\R}\Psi_{j,K}' u\,  (1-\partial_x^2)^{-1}(2u^2+u_x^2) \, dx \nonumber \\
& \le  & -\frac{c_1}{2} \int_{\R} \Psi_{j,K}'(u^2+u_x^2) +J_1+J_2+J_3 \; .
\label{go2}
\end{eqnarray}
We claim that for  $ i\in \{1,2,3\} $, it holds
\begin{equation}
J_i \le \frac{c_1}{8} \int_{\R} \Psi_{j,K}' (u^2+u_x^2) + \frac{C}{K} \, e^{-\frac{1}{K}
(\sigma_0 t+L/8)}\; . \label{go4}
\end{equation}
To handle with $ J_1 $ we   divide $ \R $ into two regions $ D_j $ and $ D_j^c $ with
 $$
 D_j=[{\tilde x_{j-1}}(t) +L/4, {\tilde x_j}(t) -L/4]
 $$
 First since from \re{eloi}, for $ x\in D_j^c $ ,
 $$
 |x-y_j(t)| \ge \frac{{\tilde x_{j}}(t)-{\tilde x_{j-1}}(t)}{2}-L/4 \ge \frac{c_j-c_{j-1}}{2} \, t +L/8\, ,
 $$
 we infer from the definition of $ \Psi $ in Section \ref{Sectmonotonie} that
 $$
 \int_{D_j^c} \Psi_{j,K}' (u^3+4 u u_x^2)  \le \frac{C}{K} \, \|u_0\|_{H^1}^3
  e^{-\frac{1}{K}(\sigma_0 t +L/8)} \; .
  $$
  On the other hand, on $ D_j $ we notice,  according to \re{distH1}, that
  \begin{eqnarray}
  \|u(t)\|_{L^\infty_{D_j}}& \le &  \sum_{i=1}^N \|\varphi_{c_i}(\cdot -{\tilde x_i(t))\|_{L^\infty}(D_j)}
   + \|u-\sum_{i=1}^N\varphi_{c_i}(\cdot -{\tilde x_i(t)})\|_{L^\infty(D_j)} \nonumber \\
    & \le & C \, e^{-L/8} +O(\sqrt{\alpha}) \; .\label{go3}
    \end{eqnarray}
Therefore, for $  \alpha $ small enough and $ L $ large enough it holds
$$
J_1 \le \frac{c_1}{8} \int_{\R} \Psi_{j,K}' (u^2+u_x^2) + \frac{C}{K} \, e^{-\frac{1}{K}(\sigma_0 t +L/8)}\; .
$$
Since $ J_2 $  can be handled in exactly the same way, it remains to treat $ J_3 $.
For this, we first notice as above that
  \begin{eqnarray}
& &\hspace*{-15mm} -\int_{D_j^c}u \Psi_{j,K}'
(1-\partial_x^2)^{-1}(2u^2+u_x^2 ) \nonumber \\
& & \le 2 \|u\|_{\infty} \sup_{x\in D_j^c}
|\Psi_{j,K}'(x-y_j(t))|\int_{\R} e^{-|x|} \ast (u^2+u_x^2 ) \, dx \nonumber \\
 &  & \le \frac{C}{K} \|u_0\|_{H^1}^3 \, e^{-\frac{1}{K}(\sigma_0 t +L/8)}\; ,
 \label{J31}
\end{eqnarray}
since
\begin{equation}
 \forall f\in L^1(\R), \quad (1-\partial_x^2)^{-1} f
  =\frac{1}{2} e^{-|x|} \ast f \; .
 \label{tytu}
 \end{equation}
Now in the region $D_j $, noticing that $ \Psi_{j,K}' $ and
$ u^2+u_x^2/2 $ are non-negative, we  get
 \begin{eqnarray}
 & & \hspace*{-15mm} -\int_{D_j} u \Psi_{j,K}'
(1-\partial_x^2)^{-1}(2u^2+u_x^2 ) \nonumber \\
 &  \le &
\|u(t)\|_{L^\infty({D_j})}\int_{D_j}\Psi_{j,K}'(
(1-\partial_x^2)^{-1}(2u^2+u_x^2) \nonumber \\
&  \le &  \|u(t)\|_{L^\infty({D_j})}\int_{\R} (2u^2+u_x^2) (1-\partial_x^2)^{-1}
\Psi_{j,K}'\; .
\end{eqnarray}
On the other hand, from the definition of $ \Psi $ in Section \ref{Sectmonotonie}
   and \re{tytu} we  infer  that for $ K\ge 4 $,
  $$
(1-\partial_x^2) \Psi_{j,K}' \ge (1-\frac{10}{K^2}) \Psi_{j,K}' \Rightarrow
(1-\partial_x^2)^{-1} \Psi_{j,K}'\le (1-\frac{10}{K^2})^{-1} \Psi_{j,K}' \; .
  $$
Therefore, taking $ K\ge 4 $ and using \re{go3}   we deduce  for $ \alpha $ small enough and $ L $ large enough that
\begin{equation}
  -\int_{D_j} u  \Psi_K'
(1-\partial_x^2)^{-1}(2u^2+u_x^2)
  \le   \frac{  c_1}{8}
\int_{\R} (u^2+u_x^2)
 \Psi_K' \; . \label{J32}
\end{equation}
This completes the proof of \re{go4}. Gathering \re{go2} and \re{go4} we  infer that
$$
\frac{d}{dt}\int_{\R} \Psi_{j,K}(u^2+u_x^2)  \, dx\le  -\frac{c_1}{8}
\int_{\R}\Psi_{j,K}' (u^2+u_x^2) +
 \frac{C}{K} \|u_0\|_{H^1}^3 \, e^{-\frac{1}{K}(\sigma_0 t +L/8)} \; .
$$
Integrating this inequality between $ 0 $ and $ t $, \re{monotonicityestim} follows.

\end{document}